\documentclass[12pt]{article}
\usepackage{latexsym}
\usepackage{amsfonts}
\usepackage{amsmath}
\usepackage{amsthm}
\usepackage{amssymb}
\usepackage{amscd}
\usepackage[dvips]{graphicx}
\usepackage{color} 
\usepackage{eepic} 
\usepackage{calrsfs}
\usepackage{txfonts}

\numberwithin{equation}{section}
\makeatletter
\newfont{\footsc}{cmcsc10 at 8truept}
\newfont{\footbf}{cmbx10 at 8truept}
\newfont{\footrm}{cmr10 at 10truept}
\newenvironment{demo}[1]{%
  \trivlist
  \item[\hskip\labelsep
        {\bf #1.}]
}{%
\hfill\qedsymbol
  \endtrivlist}

\theoremstyle{definition}
\newtheorem{theorem}{Theorem}[section]
\newtheorem{prop}[theorem]{Proposition}
\newtheorem{lemma}[theorem]{Lemma}

\def\covered{\mathinner{\mkern1mu\raise0pt\vbox{\kern7pt\hbox{$<$}}
     \mkern-4mu\raise2pt\hbox{.}\mkern2mu}}
\def\covers{\mathinner{\mkern1mu\raise0pt\vbox{\kern7pt\hbox{$>$}}
     \mkern-12mu\raise2pt\hbox{.}\mkern8mu}}

\def\defterm#1{{\sl #1}\/}

\def\operatorname#1{{\mathrm{#1}\>\!}}
\def\Bbb#1{{\mathbb{#1}}}

\def\N{\Bbb{N}}
\def\Z{\Bbb{Z}}
\def\P{{\Bbb{Z}_{>0}}}

\def\qint#1{{[{#1}]_q}}

\def\qbinom#1#2{{\left[{{#1}\atop{#2}}\right]_{q}}}

\def\rdots{\mathinner{\mkern1mu\raise0pt\vbox{\kern7pt\hbox{.}}
     \mkern2mu\raise4pt\hbox{.}\mkern2mu\raise8pt\hbox{.}\mkern1mu}}

\def\covered{\mathinner{\mkern1mu\raise0pt\vbox{\kern7pt\hbox{$<$}}
     \mkern-4mu\raise2pt\hbox{.}\mkern2mu}}
\def\covers{\mathinner{\mkern1mu\raise0pt\vbox{\kern7pt\hbox{$>$}}
     \mkern-12mu\raise2pt\hbox{.}\mkern8mu}}

\def\defterm#1{{\sl #1}\/}

\def\operatorname#1{{\mathrm{#1}\>\!}}

\def\Bbb#1{{\mathbb{#1}}}

\newcommand{\rank}{\operatorname{rank}}

\providecommand{\MR}{\relax\ifhmode\unskip\space\fi MR }

\providecommand{\href}[2]{#2}

\title{
A binomial-coefficient identity\\
 arising from \\
the middle discrete series of {${\rm SU}(2,2)$}
}

\author{
Takahiro HAYATA
\\
\small Faculty of Engineering, Yamagata University\\[-0.8ex]
\small Yonezawa, Yamagata, Japan\\[-0.8ex]
\small \texttt{hayata@yz.yamagata-u.ac.jp}
\and
Masao ISHIKAWA
\\
\small Faculty of Education, Tottori University\\[-0.8ex]
\small Koyama, Tottori, Japan\\[-0.8ex]
\small \texttt{ishikawa@fed.tottori-u.ac.jp}
}

\date{
}

\begin{document}
\maketitle

\abstract{
The aim of this paper is to give an elementary proof 
of certain identities on binomials and state an answer
to Remark~8.2 in
Takahiro Hayata, Harutaka Koseki, and Takayuki Oda, \emph{Matrix coefficients
  of the middle discrete series of {${\rm SU}(2,2)$}}, J. Funct. Anal.
  \textbf{185} (2001), 297--341.
}
\medskip

\section{Introduction}\label{sec:intro}

The aim of this paper is to show an elementary proof 
of certain identities on binomials and state an answer
to~\cite[Remark~8.2]{hko}.

The identity which we prove in this paper stems from
the representation theory of real semi-simple Lie groups, which
admits discrete series.
To describe this, let $G$ be a real semi-simple Lie group and
$K$ be its maximal compact group.
Take a unitary representation $\pi$ of $G$. Consider
the map
\[
\phi: \pi \to C^{\infty}(K\backslash G/K; \tau \otimes \tau^{*}).
\]
For a $K$-finite vector $v$, $\phi(v)$ satisfies, by definition,
$\phi(v)(kgk')=\tau(k^{-1})\tau^{*}(k')\phi(v)(g)$ $(k,k'\in K,
g\in G)$, and we say $\phi(v)$ a matrix coefficient of $\pi$ with respect to
$\tau$.
Because $G$ has the Cartan decomposition $G=KAK$ where $A$ is a
maximal split torus in $G$, the radial part,{\em i.e.}, the restriction of matrix
coefficients to $A$ is regarded as a $\tau$-valued function on
Euclidean domain.

Now we assume $\rank(G)=\rank(K)$ for $G$ to admit a discrete series representation, say $\pi$.
 If $G$ is of hermitian type and $\pi$ is
holomorphic, then $\phi(v)$ is described by the Laurent
polynomials of certain hyperbolic functions if we take its radial
part, known in the theory of Bergman kernels on the symmetric
domain.
If the Gel'fand-Kirillov dimension of $\pi$ is
enough high, the radial part can be also highly transcendental.
But the dimension is relatively low, the radial part function is expected to
be tractable. In fact, it is turned out to be feasible
when $G$ is a unitary group of degree $4$ defined by the form 
$|z_{1}|^2+|z_{2}|^2-|z_{3}|^2-|z_{4}|^2$, and $\pi$ is the second lowest
discrete series~(in~\cite{hko}, it is called \textit{a middle
  discrete series}). Because this situation is the origin of our
binomial relations, we describe details.
In this case $\dim A=2$; the radial part is a $2$-variable
function.
We find certain transform forces the function into separation of variables;
one side is described by polynomial and the other side is essentially a
Gaussian hypergeometric function ${}_{2}F_{1}$. 
The binomials $\beta_m(r,s,k,l)$
we treat here is nothing but the coefficients appearing on the
polynomial side. The desired identity rephrases 
that the value of matrix coefficients
at unit matrix is a unit.

Because we believe this kind of computation against matrix coefficients 
should work in somewhat broader contexts and likely to produce
similar identities containing involved binomial coefficients, we
hope our computation helps those who try to prove them.

\medbreak

Next, some terminology is defined before stating the main theorem.
Let $\Z$, $\N$ and $\P$ denote the set of integers, nonnegative integers
and positive integers, respectively.
If $k$ is a positive integer and $r_1$, $\dots$, $r_k$ are integers
such that $n=r_1+\cdots+r_k$ is nonnegative integer,
the \defterm{multinomial coefficient} $\binom{n}{r_1,\dots,r_k}$ is, by definition,
\begin{equation}
\begin{cases}
\frac{n!}{r_1!\cdots r_k!}
&\text{ if all $r_i\geq0$,}\\
0
&\text{ otherwise.}
\end{cases}
\end{equation}
Especially,
when $k=2$,
$\binom{n}{r}=\binom{n}{r,n-r}$ is called the \defterm{binomial coefficient}.
(For many interesting identities these famous coefficients satisfy, see \cite{K}.)
When $a$, $b$ $s$, $l$ and $m$ are integers
such that $s\geq a\geq 0$ and $s\geq l$,
we write
\begin{equation}
\beta_{m}(s,l,a,b)=\sum_{n=0}^{|b|+m}
\binom{s-a}{b_{+}+m-n}\binom{a}{b_{-}+m-n}\binom{s-l+n}{n},
\label{eq:beta}
\end{equation}
where $b_{+}$ and $b_{-}$ are defined by
\begin{align}
b_{+}+b_{-}=|b|,\qquad
b_{+}-b_{-}=b.
\end{align}
In other word
$b_{\pm}$
is defined to be
$
\frac{|b|\pm b}2.
$
The aim of this paper is to give an elementary proof of the following theorem.
\begin{theorem}
\label{th:main}
Let $s$ and $l$ be nonnegative integers such that $s\geq l\geq0$.
Let $j$ be an integer.
Then we have
\begin{align}
&\sum_{m=0}^{\lfloor (l-1)/2\rfloor}\sum_{i=2m}^{l-1}
(-2)^{i-2m}\left\{\binom{s-l}{i-2m,j-l+m,s-i-j+m} \, \beta_{m}(s,l,i+j-l,l-i)\right.
\nonumber\\
&\qquad\qquad+\left.\binom{s-l}{i-2m,j-i+m,s-l-j+m} \, \beta_{m}(s,l,l+j-i,i-l)\right\}
\nonumber\\
&\qquad\qquad+\sum_{m=0}^{\lfloor l/2\rfloor}(-2)^{l-2m}\binom{s-l}{l-2m,j-l+m,s-l-j+m} \, \beta_{m}(s,l,j,0)=\binom{s}{j}
,\label{eq:main}
\end{align}
where $\lfloor x\rfloor$ stands for the largest integer less than or equal to $x$ for any real number $x$.
Note that the left-hand side apparently include the parameter $l$,
but, eventually, it is independent of $l$ in the right-hand side.
\end{theorem}

\section{Proof of the identity}\label{sec:proof}

First we summarize certain recurrence properties of $\beta_m$ as follows.
\begin{lemma}
\label{lem:beta}
Let $s$, $l$, $a$ and $b$ be integers such that
 $s\geq l$ and $s\geq a\geq 0$.
Then the following identities hold.
\begin{enumerate}
\item[(i)]
If $b>0$, then
\begin{equation}
\label{eq01}
\beta_{m}(s,l,a,b)=\beta_{m}(s+1,l,a,b)-\beta_{m}(s+1,l,a+1,b-1).
\end{equation}
\item[(ii)]
If $a>0$ and $b\geq0$, then
\begin{equation}
\label{eq02}
\beta_{m}(s,l,a-1,b)=\beta_{m}(s+1,l,a,b)-\beta_{m-1}(s+1,l,a-1,b+1).
\end{equation}
\item[(iii)]
If $b\leq0$, then
\begin{equation}
\label{eq03}
\beta_{m}(s,l,a,b)=\beta_{m}(s+1,l,a,b)-\beta_{m-1}(s+1,l,a+1,b-1).
\end{equation}
\item[(iv)]
If $a>0$ and $b<0$, then
\begin{equation}
\label{eq04}
\beta_{m}(s,l,a-1,b)=\beta_{m}(s+1,l,a,b)-\beta_{m}(s+1,l,a-1,b+1).
\end{equation}
\end{enumerate}
\end{lemma}
\begin{demo}{Proof}
We only use the well-known recurrence equation of binomial coefficients
which reads
\begin{equation}
\binom{n+1}{r}=\binom{n}{r}+\binom{n}{r-1},
\label{eq:binom-rec}
\end{equation}
and perform direct computations to prove these identities.
First we prove (\ref{eq01}).
Applying the recurrence (\ref{eq:binom-rec}) to the first and third binomial coefficients
of (\ref{eq:beta}),
 we obtain that 
$\beta_{m}(s,l,a,b)$ equals
\begin{align*}
&\sum_{n\geq0}
\binom{s-a+1}{b+m-n}\binom{a}{m-n}\binom{s-l+n+1}{n}
-\sum_{n\geq0}
\binom{s-a}{b+m-n-1}\binom{a}{m-n}\binom{s-l+n+1}{n}\\
&\qquad\qquad-\sum_{n\geq0}
\binom{s-a}{b+m-n}\binom{a}{m-n}\binom{s-l+n}{n-1}.
\end{align*}
If we apply the recurrence $\binom{a}{m-n}=\binom{a+1}{m-n}-\binom{a}{m-n-1}$ to the second term,
then we obtain this equals
\begin{align*}
&\beta_{m}(s+1,l,a,b)
-\sum_{n\geq0}\binom{s-a}{b+m-n-1}\binom{a+1}{m-n}\binom{s-l+n+1}{n}
\\&\qquad
-\sum_{n\geq0}\binom{s-a}{b+m-n}\binom{a}{m-n}\binom{s-l+n}{n-1}
+\sum_{n\geq0}
\binom{s-a}{b+m-n-1}\binom{a}{m-n-1}\binom{s-l+n+1}{n}.
\end{align*}
The last two terms kill each other and consequently we obtain
\begin{align*}
&\beta_{m}(s,l,a,b)=\beta_{m}(s+1,l,a,b)-\beta_{m}(s+1,l,a+1,b-1).
\end{align*}
This proves the first identity.
The other identities can be proven similarly.
The details are left to the reader.
\end{demo}

Let $s$, $l$, $m$, $j$ be integers such that $s\geq l$ and $m\geq0$.
We define $\Lambda_{m}(s,l,j)$ by
\begin{align}
\Lambda_{m}(s,l,j)
&=\sum_{i=2m}^{l-1}(-2)^{i-2m}\,
\binom{s-l}{i-2m,j-l+m,s-i-j+m}\,
\beta_{m}(s,l,i+j-l,l-i)
\nonumber\\
&+\sum_{i=2m}^{l}(-2)^{i-2m}\,
\binom{s-l}{i-2m,j-i+m,s-l-j+m}\,
\beta_{m}(s,l,l+j-i,i-l).
\label{eq:Lambda}
\end{align}
Then $\Lambda_{m}(s,l,j)$ satisfies the following recurrence equation.
\begin{lemma}
\label{lem:Lambda}
Let $s$, $l$, $m$, $j$ be integers such that $s\geq l$ and $j\geq0$.
Then 
\begin{align}
&\Lambda_{m}(s,l,j)+\Lambda_{m}(s,l,j-1)=\Lambda_{m}(s+1,l,j)
+\Phi_{m}(s,l,j)-\Phi_{m-1}(s,l,j)
\label{eq:lemma}
\end{align}
where
\begin{align}
&\Phi_{m}(s,l,j)=\sum_{i=2m+1}^{l-1}(-2)^{i-2m-1}\Biggl\{\binom{s-l}{i-2m-1,j-l+m,s-i-j+m+1}\,
\beta_{m}(s+1,l,i+j-l,l-i)\nonumber\\
&\qquad\qquad+\binom{s-l}{i-2m-1,j-i+m,s-l-j+m+1}\,
\beta_{m}(s+1,l,l+j-i,i-l)\Biggr\}.
\label{eq:Phi}
\end{align}
\end{lemma}
\begin{demo}{Proof}
By (\ref{eq01}) and (\ref{eq03}),
 we obtain $\Lambda_{m}(s,l,j)$ equals
\begin{align*}
&\sum_{i=2m}^{l-1}(-2)^{i-2m} \binom{s-l}{i-2m,j-l+m,s-i-j+m}
\\&\qquad\qquad\times
\biggl\{\,\beta_{m}(s+1,l,i+j-l,l-i)
-\beta_{m}(s+1,l,i+j-l+1,l-i-1)\,\biggr\}\\
&+\sum_{i=2m}^{l}(-2)^{i-2m}
\binom{s-l}{i-2m,j-i+m,s-l-j+m}
\\&\qquad\qquad\times
\biggl\{\beta_{m}(s+1,l,l+j-i,i-l)
-\beta_{m-1}(s+1,l,l+j-i+1,i-l-1)\biggr\}.
\end{align*}
Similarly,
using (\ref{eq02}) and (\ref{eq04}),
we can rewrite $\Lambda_{m}(s,l,j-1)$ as
\begin{align*}
&\sum_{i=2m}^{l}(-2)^{i-2m}
\binom{s-l}{i-2m,j-l+m-1,s-i-j+m+1}\\
&\qquad\times 
\biggl\{\,\beta_{m}(s+1,l,i+j-l,l-i)
-\beta_{m-1}(s+1,l,i+j-l-1,l-i+1)\,\biggr\}\\
&+\sum_{i=2m}^{l-1}(-2)^{i-2m}
\binom{s-l}{i-2m,j-i+m-1,s-l-j+m+1}\\
&\qquad\times 
\biggl\{\,\beta_{m}(s+1,l,l+j-i,i-l)
-\beta_{m}(s+1,l,l+j-i-1,i-l+1)\,\biggr\}.
\end{align*}
Adding these two identities,
we obtain that
$\Lambda_{m}(s,l,j)+\Lambda_{m}(s,l,j-1)$
is equal to
\begin{align*}
&\sum_{i=2m}^{l}(-2)^{i-2m}\,
 A\,\beta_{m}(s+1,l,i+j-l,l-i)
+\sum_{i=2m}^{l-1}(-2)^{i-2m}\,
B\,\beta_{m}(s+1,l,l+j-i,i-l)\\
&-\sum_{i=2m}^{l-1}(-2)^{i-2m}
\binom{s-l}{i-2m,j-l+m,s-i-j+m}
\,\beta_{m}(s+1,l,i+j-l+1,l-i-1)\\
&-\sum_{i=2m}^{l}(-2)^{i-2m}
\binom{s-l}{i-2m,j-i+m,s-l-j+m}
\,\beta_{m-1}(s+1,l,l+j-i+1,i-l-1)\\
&-\sum_{i=2m}^{l}(-2)^{i-2m}
\binom{s-l}{i-2m,j-l+m-1,s-i-j+m+1}
\,\beta_{m-1}(s+1,l,i+j-l-1,l-i+1)\\
&-\sum_{i=2m}^{l-1}(-2)^{i-2m}
\binom{s-l}{i-2m,j-i+m-1,s-l-j+m+1}
\,\beta_{m}(s+1,l,l+j-i-1,i-l+1),
\end{align*}
where
\begin{align*}
A&=\binom{s-l}{i-2m,j-l+m,s-i-j+m}+\binom{s-l}{i-2m,j-l+m-1,s-i-j+m+1},\\
B&=\binom{s-l}{i-2m,j-i+m,s-l-j+m}+\binom{s-l}{i-2m,j-i+m-1,s-l-j+m+1}.
\end{align*}
If we replace $i$ by $i+1$ or $i-1$ in the last four terms,
then this sum becomes
\begin{align*}
&\sum_{i=2m}^{l}(-2)^{i-2m}\,
 A\,\beta_{m}(s+1,l,i+j-l,l-i)
+\sum_{i=2m}^{l-1}(-2)^{i-2m}\,
B\,\beta_{m}(s+1,l,l+j-i,i-l)\\
&-\sum_{i=2m+1}^{l}(-2)^{i-2m-1}
\binom{s-l}{i-2m-1,j-l+m,s-i-j+m+1} \, \beta_{m}(s+1,l,i+j-l,l-i)\\
&-\sum_{i=2m-1}^{l-1}(-2)^{i-2m+1}
\binom{s-l}{i-2m+1,j-i+m-1,s-l-j+m} \, \beta_{m-1}(s+1,l,l+j-i,i-l)\\
&-\sum_{i=2m-1}^{l-1}(-2)^{i-2m+1}
\binom{s-l}{i-2m+1,j-l+m-1,s-i-j+m} \, \beta_{m-1}(s+1,l,i+j-l,l-i)\\
&-\sum_{i=2m+1}^{l}(-2)^{i-2m-1}
\binom{s-l}{i-2m-1,j-i+m,s-l-j+m+1} \, \beta_{m}(s+1,l,l+j-i,i-l).
\end{align*}
Using
\begin{align*}
&A+\binom{s-l}{i-2m-1,j-l+m,s-i-j+m+1}=\binom{s-l+1}{i-2m,j-l+m,s-i-j+m+1},\\
&B+\binom{s-l}{i-2m-1,j-i+m,s-l-j+m+1}=\binom{s-l+1}{i-2m,j-i+m,s-l-j+m+1},
\end{align*}
we see that 
$\Lambda_{m}(s,l,j)+\Lambda_{m}(s,l,j-1)$
is equal to
\begin{align*}
&\sum_{i=2m}^{l-1}(-2)^{i-2m}
\binom{s-l+1}{i-2m,j-l+m,s-i-j+m+1} \,
\beta_{m}(s+1,l,i+j-l,l-i)\\
&+\sum_{i=2m}^{l}(-2)^{i-2m}
\binom{s-l+1}{i-2m,j-i+m,s-l-j+m+1} \,
\beta_{m}(s+1,l,l+j-i,i-l)\\
&+\sum_{i=2m+1}^{l-1}(-2)^{i-2m-1}
 \binom{s-l}{i-2m-1,j-l+m,s-i-j+m+1} \,
\beta_{m}(s+1,l,i+j-l,l-i)\\
&-\sum_{i=2m-1}^{l-1}(-2)^{i-2m+1}
 \binom{s-l}{i-2m+1,j-i+m-1,s-l-j+m} \,
\beta_{m-1}(s+1,l,l+j-i,i-l)\\
&-\sum_{i=2m-1}^{l-1}(-2)^{i-2m+1}
 \binom{s-l}{i-2m+1,j-l+m-1,s-i-j+m} \,
\beta_{m-1}(s+1,l,i+j-l,l-i) \\
&+\sum_{i=2m+1}^{l-1}(-2)^{i-2m-1}
 \binom{s-l}{i-2m-1,j-i+m,s-l-j+m+1} \,
 \beta_{m}(s+1,l,l+j-i,i-l),
\end{align*}
which is equal to the right-hand side of \eqref{eq:lemma}.
This completes the proof of the lemma.
\end{demo}

\bigbreak
\noindent
\begin{demo}{Proof of Theorem~\ref{th:main}}
Assume $s\geq l\geq0$.
If we put
\begin{equation*}
\Gamma(s,l,j)=\sum_{m=0}^{\infty}\Lambda_{m}(s,l,j),
\end{equation*}
then,
by \eqref{eq:lemma},
it is easy to see that
\begin{equation}
\Gamma(s+1,l,j)=\Gamma(s,l,j)+\Gamma(s,l,j-1)
\label{eq:rec-eq}
\end{equation}
holds.
In addition,
if $j<0$ or $j>s$,
then we have $\Lambda_{m}(s,l,j)=0$ for all $m\geq0$
since the multinomial coefficients vanish
in the definition \eqref{eq:Lambda}.
If $j=0$,
then we also have
\[
\Lambda_{m}(s,l,j)
=\begin{cases}
\beta_{0}(s,l,l,-l)=1
&\text{ if $m=0$,}\\
0
&\text{ if $m>0$.}
\end{cases}
\]
Hence, we have
\begin{equation}
\Gamma(s,l,j)
=\begin{cases}
0
&\text{ if $j<0$ or $j>s$,}\\
1
&\text{ if $j=0$.}
\end{cases}
\label{eq:initial}
\end{equation}
From \eqref{eq:rec-eq} and \eqref{eq:initial},
we conclude that $\Gamma(s,l,j)=\binom{s}{j}$.
This completes the proof.
\end{demo}

\section{Concluding remarks}

An interesting question we can ask is 
``Can one make a $q$-analogue of the identity \eqref{eq:main}?''.
(For $q$-series, the reader can refer to \cite{A1}.)
We had a trial in this direction which is not yet successful.
For example,
define $\beta_{m}(s,l,a,b)$ by
\begin{equation}
\beta_{m}(s,l,a,b)=\sum_{n=0}^{|b|+m}q^{n(n-|b|+l-2m)}
\qbinom{s-a}{b_{+}+m-n}\qbinom{a}{b_{-}+m-n}\qbinom{s-l+n}{n},
\label{eq:beta-q}
\end{equation}
where
\begin{equation*}
\qbinom{r_{1}+\cdots+r_{k}}{r_{1},\dots,r_{k}}
=\begin{cases}
\frac{\qint{r_{1}+\cdots+r_{k}}!}{\qint{r_1}!\cdots\qint{r_k}!}
&\text{ if all $r_i\geq0$,}\\
0
&\text{ otherwise,}
\end{cases},
\qquad\qquad
\qbinom{n}{r}=\qbinom{n}{r,n-r},
\end{equation*}
with $\qint{n}!=(1+q)\cdots(1+q+\cdots+q^{n-1})$.
Then
one can prove that
$\beta_{m}(s,l,a,b)$ satisfies the following simple recurrence equations.


\begin{prop}
\label{lem:beta-q}
Let $s$, $l$, $a$ and $b$ be integers such that
 $s\geq l$ and $s\geq a\geq 0$.
Then the following identities hold.
\begin{enumerate}
\item[(i)]
If $b>0$, then
\begin{equation}
\label{eq01q}
\beta_{m}(s,l,a,b)=\beta_{m}(s+1,l,a,b)-q^{s-a-b-m+1}\beta_{m}(s+1,l,a+1,b-1)
\end{equation}
\item[(ii)]
If $a>0$ and $b\geq0$, then
\begin{equation}
\label{eq02q}
\beta_{m}(s,l,a-1,b)=\beta_{m}(s+1,l,a,b)-q^{a-m}\beta_{m-1}(s+1,l,a-1,b+1)
\end{equation}
\item[(iii)]
If $b\leq0$, then
\begin{equation}
\label{eq03q}
\beta_{m}(s,l,a,b)=\beta_{m}(s+1,l,a,b)-q^{s-a-m+1}\beta_{m-1}(s+1,l,a+1,b-1)
\end{equation}
\item[(iv)]
If $a>0$ and $b<0$, then
\begin{equation}
\label{eq04q}
\beta_{m}(s,l,a-1,b)=\beta_{m}(s+1,l,a,b)-q^{a+b-m}\beta_{m}(s+1,l,a-1,b+1)
\end{equation}
\end{enumerate}
\end{prop}
Nevertheless,
at this point,
we don't know how to define
$\Lambda_{m}(s,l,j)$ which would have a simple recurrence equation.

\medbreak

Another interesting question we can ask is the following.
One can see that the left-hand side of \eqref{eq:main} is a sum
(a double sum or triple sum),
but the right-hand side is so simple, i.e.,
just a binomial coefficient $\binom{s}{j}$.
So one may ask whether any of the algorithms
such as the WZ algorithm  (see \cite{Z2})
would be able to handle it?
(We would like to thank to Prof. C.~Krattenthaler for his helpful comment).

\end{document}